\documentclass[preprint,10pt]{elsarticle}
\usepackage{amssymb}
\usepackage{amsmath}
\usepackage{amsthm}
\usepackage{amsfonts}
\usepackage[all]{xy}
\usepackage{color}
\usepackage{filecontents}
\usepackage{natbib}
\usepackage{tikz}
\usetikzlibrary{positioning,automata}
\usepackage{lineno}
\usepackage{ mathrsfs }
\usepackage{bm}

\theoremstyle{definition}

\newtheorem{theorem}{Theorem}[section]

\newtheorem{corollary}{Corollary}[section]

\newtheorem{definition}{Definition}[section]
\newtheorem{example}{Example}[section]

\DeclareMathOperator{\supp}{supp}

\newcommand{\R}{\mathbb{R}}  % The real numbers.
\newcommand{\N}{\mathbb{N}_0}  % The Natural numbers.
\numberwithin{equation}{section}

\begin{document}

\begin{frontmatter}

\title{Catalan-like numbers and Hausdorff moment sequences }

\author[label1]{Hayoung Choi}
\ead{hchoi@shanghaitech.edu.cn}
\address[label1]{School of Information Science and Technology, ShanghaiTech University, Pudong district, Shanghai 200031, China}

\author{Yeong-Nan Yeh\fnref{label2}}
\ead{mayeh@math.sinica.edu.tw}
\address[label2]{Institute of Mathematics, Academia Sinica, Taipei, Taiwan}

\author{Seonguk Yoo\fnref{label3}\corref{corr}}
\ead{seyoo@gnu.ac.kr}
\address[label3]{Department of Mathematics Education and RINS, Gyeongsang National University, Jinju, Republic of Korea}

\cortext[corr]{Corresponding author}

\begin{abstract}
In this paper we show that many well-known counting coefficients, including the Catalan numbers, the Motzkin numbers, the central binomial coefficients, the central Delannoy numbers are Hausdorff moment sequences in a unified approach.
In particular we answer a conjecture of Liang at al. which such numbers have unique representing measures.
The smallest interval including the support of representing measure is explicitly found. 
Subsequences of Catalan-like numbers are also considered. 
We provide a necessary and sufficient condition for a pattern of subsequences that if sequences are the Stieltjes Catalan-like numbers, then their subsequences are Stieltjes Catalan-like numbers.
Moreover, a representing measure of a linear combination of consecutive Catalan-like numbers is studied. 

%We show that well-known Catalan-like numbers are Hausdorff moment sequences, which are determinate.

\end{abstract}

\end{frontmatter}

%\linenumbers

%% Introduction
\section{Introduction}

Let $\N$ (resp., $\R$) be the set of nonnegative integers (resp., real numbers). %and $\R_+$ be the set of nonnegative real numbers.
%Let us denote by $\N^n$ the set of all multi-indices  $\bm{\alpha} := (\alpha_1, \ldots, \alpha_n)$, i.e., $\alpha_j \in \N$ for all $j=1,\ldots,n$.
%A \emph{sequence} is a function from $\N^n$ to $\R$, denoted by $y:=(y_{\bm{\alpha}})_{\bm{\alpha} \in \N^n}$.
Let $K$ be a closed subset of $\R$. A nonnegative Borel measure $\mu$ on $\R$ is called a \emph{$K$-measure} if its support, denoted by $\supp(\mu$), is contained in $K$.
The symbol $\R[x]$ denotes the ring of polynomials in $x$ with real coefficients.
%For $\bm{x} \in \R^n$ and $\bm{\alpha} \in \N^n$, denote $\bm{x}^{\bm{\alpha}}=x_1^{\alpha_1} \cdots x_n^{\alpha_n}$. 
The integral $\int_{K}x^{n} d\mu$, if it exists, is called the \emph{$n$-th moment of the measure $\mu$}. 
A sequence $y=(y_n)_{n \geq 0}$ is said to \emph{admit a $K$-measure $\mu$} if 
\begin{equation}\label{eq:moments}
y_{n} = \int_{K} x^{n} d\mu, \quad \text{for all } 
n \in \N.
\end{equation}
Such $\mu$ is called a \emph{ $K$-representing measure for $y$} and $y$ is called a \emph{$K$-moment sequence}. When $K=\R$ (resp. $K=[0,\infty)$, $K=[a,b]$), the sequence $y$ is also called a \emph{Hamburger(resp. Stieltjes, Hausdorff) moment sequence}. A moment sequence is called \emph{determinate}, if there is a unique representing measure  such that \eqref{eq:moments} holds; otherwise it is called \emph{indeterminate}. 
For more information, see references \cite{book:momentproblem,book:Shohat,book:Akhiezer,Curto97} and references therein.

Given sequence $y=(y_n)_{n \geq 0}$, we denote 
\begin{equation}%\label{Hankelmatrix}
H_m(y):=
\begin{bmatrix}
y_{0} & y_{1} & \cdots & y_{m} \\
y_{1} & y_{2} & \cdots & y_{m+1} \\
\vdots & \vdots & \ddots & \vdots \\
y_{m} & y_{m+1} & \cdots & y_{2m} \\
\end{bmatrix}, \quad
\widetilde{H}_m(y):=
\begin{bmatrix}
y_{1} & y_{2} & \cdots & y_{m+1} \\
y_{2} & y_{3} & \cdots & y_{m+2} \\
\vdots & \vdots & \ddots & \vdots \\
y_{m+1} & y_{m+2} & \cdots & y_{2m+1} \\
\end{bmatrix}.
\end{equation}
Denote $\Delta_m(y) := \det(H_m(y))$ for all $m\in \N$.
It is well known that $y$ is a Hamburger moment sequence if and only if $H_m(y)$ is positive semidefinite for all $m \in \N$, and
$y$ is a Stieltjes moment sequence if and only if both $H_m(y)$ and $\widetilde{H}_m(y)$ are positive semidefinite for all $m \in \N$, equivalently, $H_m(y)$ is totally positive for all  $m \in \N$.

% Catlan-like numbers
Aigner\cite{AIGNER199933,Aigner2001} introduced a unified approach for many well-known combinatorial sequences including the Catalan numbers, the Ridordan numbers, the Fine numbers, the Motzkin numbers, the Schr\'{o}der numbers, and so on.

Let $\sigma=(s_k)_{k\geq 0}$ and $\tau = (t_k)_{k\geq 1}$ be two sequences of real numbers with $t_{k+1} \neq 0$ for all $k\in \N$ and define an infinite lower triangle matrix 
$R:= R^{\sigma,\tau} = [r_{n,k}]_{n,k \geq 0}$ by the recurrence relations
\begin{equation}\label{recursive_matrix}
r_{0,0}= 1, \quad	 r_{n+1,k} = r_{n,k-1} + s_k r_{n,k} + t_{k+1} r_{n,k+1},
\end{equation}
where 
$r_{n,k}=0$ unless $n \geq k \geq 0$. Clearly, all $r_{n,n} =1$. 
$R$ is called the {\emph{recursive matrix}} and $r_n = r_{n,0}$ is the {\emph{$n$th Catalan-like numbers }}corresponding to $(\sigma,\tau)$.

\begin{example}\label{ex:Catalan-like}
The following well-known counting coefficients are Catalan-like numbers.
\begin{itemize}
\item[$(i)$] the Catalan numbers $C_n$ when $\sigma  = (1,2,2,\ldots)$ and $\tau  = (1,1,1,\ldots)$;
\item[$(ii)$] the shifted Catalan numbers $C_{n+1}$ when $\sigma = (2,2,2,\ldots)$ and $\tau = (1,1,1,\ldots)$;
\item[$(iii)$] the Motzkin numbers $M_n$ when $\sigma = \tau = (1,1,1, \ldots)$;
\item[$(iv)$] the central binomial coefficients $\binom{2n}{n}$ when $\sigma=(2,2,2,\ldots)$ and $t=(2,1,1,\ldots)$;
\item[$(v)$] the central trinomial coefficients $T_n$ when $\sigma = (1,1,1, \ldots)$ and $t = (2,1,1,\ldots)$;
\item[$(vi)$] the central Delannoy numbers $D_n$ when $\sigma = (3,3,3, \ldots)$ and $\tau = (4,2,2, \ldots)$;
\item[$(vii)$] the large Schr\'{o}der numbers $r_n$ when $\sigma = (2,3,3,\ldots)$ and $\tau = (2,2,\ldots)$;
\item[$(viii)$] the little Schr\'{o}der numbers $S_n$ when $\sigma = (1,3,3,\ldots)$ and $\tau = (2,2,\ldots)$;
\item[$(ix)$] the Fine numbers $F_n$ when $\sigma = (0,2,2,\ldots)$ and $t = (1,1,1,\ldots)$;
\item[$(x)$] the Riordan numbers $R_n$ when $\sigma = (0,1,1,\ldots)$ and $t = (1,1,1,\ldots)$.
\item[$(xi)$] the (restricted) hexagonal numbers $h_n$ when $\sigma = (3, 3, 3, \ldots)$ and $\tau = (1, 1, 1, \ldots)$;
\end{itemize}
\end{example}

Liang at al.\cite{LIANG2016484} showed that these types of sequences are Stieltjes moment sequences 
in unified setting.
Wang\cite{WANG2016115} showed that Stieltjes moment sequences are infinitely log-convex.
Chen at al. \cite{CHEN2015383} presented some sufficient conditions such that the recursive matrix is totally positive. They also proved that many well-known sequences are log-convex\cite{CHEN201568}.
For more recent work, see \cite{WANG20181264,Liang2018}.

\cite{LIANG2016484} remarked that it is questionable if many well-know sequences are determinate. 
To check whether a given moment sequence is determinate is an important problem, but very difficult.
The reader is referred to \cite[Chapter 4]{book:momentproblem} for a comprehensive study and \cite{Lin2017} for recent results.

It is well-known that Hausdorff moment sequences are determinate.
Answering this conjecture, our main result is the following.
 \begin{theorem}
 The Catalan-like numbers corresponding to $ \sigma = (p, s, s, \ldots)$ and $\tau = (q, t, t, \ldots)$ with $p>s-2\sqrt{t}$ and $\max\{q,t\} < s+2\sqrt{t}$ are
 $[s-2\sqrt{t}, s+2\sqrt{t}]$-moment sequences with $\big| \supp(\mu) \big| = \infty$.
 \end{theorem}
Many well-known Catalan-like numbers have the form $ \sigma = (p, s, s, \ldots)$ and $\tau = (q, t, t, \ldots)$
It shows that all sequences in Example \ref{ex:Catalan-like} are Hausdorff moment sequences, implying that each has the unique representing measure.
Specifically, the following are shown.
\begin{itemize}
\item[$(1)$] The Catalan numbers $C_n$,  the shifted Catalan numbers $C_{n+1}$,
the central binomial coefficients $\binom{2n}{n}$, and the Fine numbers $F_n$ 
are $[0,4]$-moment sequences.

\item[$(2)$] The Motzkin numbers $M_n$,
 the central trinomial coefficients $T_n$, and
the Riordan numbers $R_n$ are $[-1,3]$-moment sequences.

\item[$(3)$] The central Delannoy numbers $D_n$,
the large Schr\'{o}der numbers $r_n$, and
the little Schr\'{o}der numbers $S_n$
are $[3-2\sqrt{2},3+2\sqrt{2}]$-moment sequences.

 \item[$(4)$] The (restricted) hexagonal numbers $h_n$ are $[1,5]$-moment sequences.
\end{itemize}

Using these results, we also provide a necessary condition for each Catalan-like numbers to be a Hausdorff moment sequence 
in Corollary \ref{cor:necessary_condition}. For such necessary condition it is required to find the closed interval including the support.

Finding integral representations of well-known counting numbers has been studied(\cite{math5030040} and reference therein).
Finding a correct interval of the support of the representing measure will help one find an integral representations.

Next new sequences from Catalan-like numbers are considered. We mainly consider two different type of new sequences: (1) subsequences of Catalan-like numbers; (2) linear combinations of Catalan-like numbers.

It is shown that a necessary and sufficient condition for $n_k$ which the subsequence $(y_{n_k})_{n\geq 0}$ is Stieltjes moment sequences for all Stieltjes moment sequences $(y_n)_{n\geq 0}$ is
$$n_k=dk+\ell \quad \text{for all } k \in \N,$$
where $d,\ell\in \N$.

Bouras\cite{Bouras2013} considered the determinant of $H_m(z)$ where $z=(z_n)_{n\geq 0}$ is defined as a linear combination of three successive shifted Catalan numbers
$$
z_n = \alpha_0 C_{n+k} + \alpha_1 C_{n+k+1} + \alpha_2 C_{n+k+2} \quad \text{for all }n\in \N,
$$
where $\alpha_0, \alpha_1, \alpha_2, \in \R$ and $k$ is an arbitrary positive integer.
It was shown that such linear combination can be expressed in terms of the moments of a linear functional related to Jacobi linear functional, based on 
the well-known relation between orthogonal polynomials and Hankel determinants.

Mu, Wang, and Yeh (2017) considered this in more general setting \cite{MU20173097}. They unify many known results of Hankel determinant evaluations for well-known counting numbers. They show an explicit form for the determinant of a linear combination of consecutive Catalan-like numbers.
We focus on the representing measures of the new sequences instead of the determinants.
For Catalan-like numbers $r=(r_n)_{n \geq 0}$, we consider a new sequence $\tilde{r} = (\tilde{r}_n)_{n \geq 0}$ defined by
$$
 \tilde{r}_n= \alpha_0 r_{dn+\ell} +  \alpha_1 r_{dn+\ell+1} +  \alpha_2 r_{dn+\ell+2} + \cdots + +  \alpha_m r_{dn+\ell+m},
 $$
where $\alpha_i \in \R$ is given for all $0\leq i \leq m$ and $d,\ell \in \N$.
%\begin{theorem}
%Suppose that  $\alpha_i$ satisfies that 
%$$\alpha_0 +  \alpha_1 t +  \cdots + \alpha_m t^m  \geq 0 \quad \text{for all } a \leq t \leq b,$$ If $r$ is a $[a,b]$-moment sequence, then
%$\tilde{r}$ is a $[a,b]$-moment sequence.
%\end{theorem}

Finding the the smallest closed interval including the support of representing measure for each Catalan-like numbers is required to study this type of new sequences.

\section{Preliminaries}

In this section we introduce powerful tools for a study of one-dimensional moment problems: Riesz functional and orthogonal polynomials.
Readers are referred to the reference  \cite{Chihara78,book:momentproblem} for deeper treatment of the results.

For a sequence, $y=(y_n)_{n \geq 0}$, define a \emph{Riesz functional} $\mathcal{L}_y$ acting on $\R[x]$ as 
\begin{equation*}
\mathcal{L}_y\Big[ \sum c_n x^n \Big] =  \sum c_n y_n.
\end{equation*}
We write simply $\mathcal{L}$ instead of $\mathcal{L}_y$ when it is understood well.

%For a sequence, $y=(y_n)_{n\geq 0}$, define a \emph{Riesz functional} $L_y$ acting on $\R[x]$  as 
%\begin{equation}
%L_y\bigg( \sum  p_n x^n \bigg) :=  \sum p_n y_n.
%\end{equation}
We say that $\mathcal{L}$ is \emph{$K$-positive} if 
\begin{equation}
\mathcal{L}_y (p) \geq 0 \quad \forall p\in \R[x]: p|_K \geq 0.
\end{equation}
When $K=\R$, we call positive instead of $K$-positive.

The $K$-positivity of $\mathcal{L}_y$ is a necessary condition for $y$ to admit a $K$-measure.
Conversely, the classical theorem of M. Riesz\cite{Riesz23} provides a fundamental existence criterion for $K$-representing measures and Haviland\cite{Haviland36} provides the generalization in $\R^n$.
\begin{theorem}[Riesz-Haviland]\label{thm:RH_full}
A sequence $y=(y_n)_{n \geq 0}$ admits a representing measure supported in the closed set $K\subset \R$ if and only if $\mathcal{L}_{y}$ is $K$-positive.
\end{theorem}

From now on, we will collect well-known results about Riesz functional and orthogonal polynomials.

%\begin{definition}
A sequence $\{ P_n(x)\}_{n\geq 0}$ is called an {\emph{orthogonal polynomial sequence}}(in short, OPS) with respect to $\mathcal{L}$ if it satisfies that 
\begin{equation*}
\deg{(P_n)} = n \quad \text{and}\quad \mathcal{L}[P_m P_n]= K_n\delta_{mn}~ (K_n \neq 0) \quad \text{for all }m,n\in \N.
\end{equation*}
When $K_n =1$ for all $n\in\N$, such an OPS is called an {\emph{orthonormal polynomial sequence}}.
%\end{definition}

%\begin{theorem}\label{thm:quasi-definite}
%For a sequence $y=(y_n)_{n\geq 0}$, let $\mathcal{L}_y$ be a Riesz functional. A necessary and sufficient condition for the existence of an OPS for $\mathcal{L}_y$ is 
%$$\Delta_n(y) \neq 0 \quad \text{for all } n\in \N.$$ 
%\end{theorem}
%Such $\mathcal{L}_y$ is called {\emph{quasi-definite}}.
%
%\begin{theorem}\label{thm:coefficient}
%Let $\{ P_n(x)\}_{n\geq 0}$ be an OPS with respect to quasi-definite $\mathcal{L}$. Then every polynomial $\pi(x)$ of degree $n$ can be written as 
%\begin{equation*}
%\pi(x) = \sum_{k=0}^n c_k P_k(x),
%\end{equation*}
%where $c_k = \dfrac{\mathcal{L}[\pi(x) P_k(x)]}{\mathcal{L}[P_k^2(x)]}$ for all $k=0,1,\cdots, n$.
%\end{theorem}
%
There exists an explicit formulas for the orthogonal polynomial sequences (see \cite[Proposition 5.3]{book:momentproblem}).
\begin{theorem}\label{thm:OPS_form}
For a sequence $y=(y_n)_{n\geq 0}$, let $\mathcal{L}_y$ be a Riesz functional. Then the monic OPS for $\mathcal{L}_y$ is 
expressed as 
\begin{equation*}%\label{Hankelmatrix}
P_n(x)=
\frac{1}{\Delta_{n-1}(y)}
\det
\begin{bmatrix}
y_{0} & y_{1} & \cdots & y_{n} \\
y_{1} & y_{2} & \cdots & y_{n+1} \\
\vdots & \vdots & \ddots & \vdots \\
y_{n-1} & y_{n} & \cdots & y_{2n-1} \\
1 & x & \cdots & x^n \\
\end{bmatrix},
\end{equation*}
provided that $\Delta_n(y) \neq 0$ for all $n\in \N.$.
\end{theorem}
In fact, the condition $\Delta_n(y) \neq 0$ for all $n\in \N.$ is a necessary and sufficient condition for the existence of an OPS for $\mathcal{L}_y$. 
Such $\mathcal{L}_y$ is called {\emph{quasi-definite}}.

%See \cite[Chapter 1]{book:Akhiezer} or \cite[Exercise 3.1, Chapter 1]{Chihara78}

%\begin{theorem}
%Let $\{ P_n(x)\}_{n\geq 0}$ be a monic OPS with respect to $\mathcal{L}_y$.
%Then, 
%\begin{equation}
%\mathcal{L}[P_n^2(x)] = \frac{\Delta_n(y)}{\Delta_{n-1}(y)} \quad \text{for all } n\in \N.
%\end{equation}
%\end{theorem}
%\begin{proof}
%By Theorem \ref{thm:coefficient}, we have 
%\begin{equation*}
%\mathcal{L}[P_n^2(x)] = \mathcal{L}[P_n(x) x^n] =   \frac{\Delta_n(y)}{\Delta_{n-1}(y)}.
%\end{equation*}
%\end{proof}

\begin{theorem}[\cite{Chihara78}, Theorem 4.1, Chapter 1]\label{thm:recurrence_OPS}
Let $\mathcal{L}_y$ be a quasi-definite Riesz functional and let $\{ P_n(x)\}_{n\geq 0}$ be the corresponding monic OPS with respect to $\mathcal{L}_y$. Then there exist $\sigma=(s_k)_{k\geq 0}$ and $\tau = (t_k)_{k\geq 1}$ with $t_{k+1} \neq 0$ for all $k\in \N$ such that 
\begin{equation}\label{eq:recurrence_OPS}
P_{k+1}(x) = (x-s_k)P_k(x) - t_k P_{k-1}(x) \quad \text{for all }k\in \N,
\end{equation}
where we definite $P_{-1}(x)=0$ and $t_0$ is arbitrary.
Furthermore, for each $k\in \N$
\begin{equation}
s_k = \frac{\mathcal{L}_y[xP_k^2(x)]}{\mathcal{L}_y[P_k^2(x)]} \quad \text{and} \quad t_{k+1} = \frac{\mathcal{L}_y[P_{k+1}^2(x)]}{\mathcal{L}_y[P_{k}^2(x)]} = \frac{\Delta_{k-1}(y) \Delta_{k+1}(y)}{(\Delta_{k}(y))^2}.
\end{equation}
Moreover, if $\mathcal{L}_y$ is positive-definite (i.e., $\mathcal{L}_y[\pi(x)]>0$ for all nonzero positive polynomial $\pi(x)$), then $s_k \in \R$ and $t_{k+1} > 0$ for all $k\in \N$. 
\end{theorem}

\begin{theorem}\label{thm:Favard}
{\bf{(Favard's Theorem)}}
Let $\sigma=(s_k)_{k\geq 0}$ and $\tau = (t_k)_{k\geq 1}$ be arbitrary sequences of complex numbers and let $\{ P_n(x)\}_{n\geq 0}$ be defined by the recurrence formula
\begin{equation}
P_{-1}(x)=0,\quad P_0=1,\quad P_{k+1}(x) = (x-s_k)P_k(x) - t_k P_{k-1}(x) \quad \text{for all }k\in \N,
\end{equation}
Then, there exists an unique Riesz functional $\mathcal{L}$ such that
\begin{equation}\label{eq:Riesz_condition}
\mathcal{L}[1]=t_0,\quad \mathcal{L}[P_m(x) P_n(x)] = 0 \quad \text{for all } m,n \in \N \text{ with }m \neq n.
\end{equation} 
Moreover, $\mathcal{L}$ is quasi-definite and $\{ P_n(x)\}_{n\geq 0}$ is the corresponding monic OPS if and only if $t_{k+1} \neq 0$ for all $k\in \N$. 
$\mathcal{L}$ is positive-definite if and only if $s_k\in \R$ and $t_{k+1}>0$ for all $k\in \N$.
\end{theorem}

%Given $\tau = (t_k)_{k\geq 1}$, we set $\hat{t}_n = t_1 t_2 \ldots t_n$ $(n\geq 1)$, $\hat{t}_0=1$.
%For $R:= R^{\sigma,\tau} = [r_{n,k}]_{n,k \geq 0}$, the submatrix of $R$ consisting of rows and columns $0$ to $m$ is denoted by $R_m$.
%\begin{lemma}\label{lemma:fundamental1}
%For $R:= R^{\sigma,\tau} = [r_{n,k}]_{n,k \geq 0}$, we have 
%\begin{equation*}
%R_m T_{m} R_m^\top = H_{m}(r),
%\end{equation*}  
%where $r=(r_n)_{n\geq 0}$ and
%\begin{equation*}
%T_m= 
%\begin{bmatrix}
%\hat{t}_0 &	 & 	& 	& \\
%       & \hat{t}_1 &     & 	\text{\huge $0$}  &\\
%       &	& \hat{t}_2 & 	&\\
%       & \text{\huge $0$}	& 	& \ddots &\\
%       &	& 	&  & \hat{t}_m
%\end{bmatrix}.
%\end{equation*}
%\end{lemma}
%\begin{proof}
%Note  \cite{Aigner2001} that for all $m,n\in \N$, 
%\begin{equation*}
%\sum_k r_{m,k} r_{n,k} \hat{t}_k = r_{m+n,0}=r_{m+n}.
%\end{equation*}
%\end{proof}

%
%\begin{corollary}
%$R:= R^{\sigma,\tau} = [r_{n,k}]_{n,k \geq 0}$ is a recursive matrix if and only if $RT_{\infty}R^\top = H_\infty$ with $\hat{t}_n \neq 0$ for all $n$, $\hat{t}_0=1$. The sequences $\sigma$ and $\tau$ are given by $s_k = r_{k+1,k} - r_{k,k-1}$ and $t_k = T_k/T_{k-1}$ for all $k\in \N$.
%\end{corollary}

%\begin{corollary}\label{cor:catalan-like-det}
%A sequence  $y=(y_n)_{n\geq 0}$ is Catalan-like if and only if $\Delta_m(y) \neq 0$ for all $m\in \N$.
%\end{corollary}
%\begin{proof}
%Since $\det{(R_m)}=1$, by Lemma \ref{lemma:fundamental1}, it holds that $\det{(T_m)} = \Delta_m(y)$. 
%\end{proof}

Aigner\cite{Aigner2001} showed an interesting connection between recursive matrices and coefficients of monic OPS.
%For $R= R^{\sigma,\tau}$, let $U=R^{-1}= [u_{n,k}]_{n,k \geq 0}$ be the inverse matrix of $R$. Clearly, $U$ is also low triangular with diagonal equal to 1, and satisfies the recurrence relations
%\begin{equation*}
%u_{0,0}= 1,\quad  u_{n+1,k} = u_{n,k-1} - s_n u_{n,k} - t_{n} u_{n-1,k} \quad \text{for all } n\in \N.
%\end{equation*}
%Define the polynomial sequence $\{ P_n(x)\}_{n\geq 0}$ by
%\begin{equation*}
%P_n(x) = \sum_{k=0}^n u_{n,k} x^k.
%\end{equation*}
%Then it follows that 
%\begin{equation}\label{eq:poly_recurrence}
%P_{-1}(x)=0,\quad P_0=1,\quad P_{n+1}(x) = (x-s_n)P_n(x) - t_n P_{n-1}(x) \quad \text{for all }n\in \N.
%\end{equation}
%Hence we obtain the fundamental equivalence between $R^{\sigma,\tau}$ and  $\{ P_n(x)\}_{n\geq 0}$ as follows:
%\begin{equation*}
% \{ P_n(x)\}_{n\geq 0} \text{ satisfies } (\ref{eq:poly_recurrence}) \Longleftrightarrow UHU^\top = T   \Longleftrightarrow RTR^\top = H \Longleftrightarrow R^{\sigma,\tau}  \text{ is a recursive matrix}
%\end{equation*}

\begin{theorem}\label{thm:main1}
Let $y=(y_n)_{n\geq 0}$ be a sequence. Then the following are equivalent:
\begin{enumerate}
\item[$(1)$]  $y$ is the Catalan-like numbers corresponding to $(\sigma,\tau)$;
\item[$(2)$]  $\Delta_m (y) \neq 0$ for all $m\in \N$ ;
\item[$(3)$] There exists a recursive matrix $R^{\sigma,\tau} = [r_{n,k}]_{n,k \geq 0}$ such that $r_{n,0} =y_n$ for all $n\in \N$;
\item[$(4)$] There exists an monic OPS, $\{ P_n(x)\}_{n\geq 0}$, with respect to $\mathcal{L}_y$; 
\item[$(5)$] There exist constants $\sigma=(s_k)_{k\geq 0}$ and $\tau = (t_k)_{k\geq 1}$ with $t_{k} \neq 0$  satisfying \eqref{eq:recurrence_OPS};
%\item[$(6)$] There exists constants $\sigma=(s_k)_{k\geq 0}$ and $\tau = (t_k)_{k\geq 1}$ satisfying \eqref{recursive_matrix};
\item[$(6)$] There exists a quasi-definite Riesz functional $\mathcal{L}_y$. 
%satisfying \eqref{eq:Riesz_condition}.
\end{enumerate}
\end{theorem}
\begin{proof}
Let
$\tau = (t_k)_{k\geq 1}$ be two sequences of real numbers with $t_{k+1} \neq 0$ for all $k\in \N$.
By  the observation in \cite{Aigner2001} (1) and (2) are equivalent.
By the definition of the Catalan-like numbers, (1) and (3) are equivalent.
By Theorem 3.1, Chapter 1\cite{Chihara78}, (2) and (4) are equivalent.
By Theorem \ref{thm:recurrence_OPS} and \ref{thm:Favard}, (4) and (5) are equivalent.
By the definition of quasi-definite Riesz functional, (2) and (6) are equivalent.
\end{proof}

\begin{theorem}\label{thm:main2}
Let $y=(y_n)_{n\geq 0}$ be a sequence. Then the following are equivalent:
\begin{enumerate}
\item[$(1)$]  $y$ is a positive sequence and the Catalan-like numbers;
\item[$(2)$]  $y$ is a positive definite sequence;
\item[$(3)$]  $\Delta_m (y) > 0$ for all $m\in \N$;
\item[$(4)$] $y$ admits a $K$-measure $\mu$ with $K \subseteq \R$ such that $\big| \supp(\mu) \big| = \infty$;
\item[$(5)$] There exists a recursive matrix $R^{\sigma,\tau} = [r_{n,k}]_{n,k \geq 0}$ such that $r_{n,0} =y_n$ and $t_k >0$ for all $n,k$;
\item[$(6)$] There exists an monic OPS, $\{ P_n(x)\}_{n\geq 0}$, with respect to positive-definite $\mathcal{L}_y$;
\item[$(7)$] There exist constants $\sigma=(s_k)_{k\geq 0}$ and $\tau = (t_k)_{k\geq 1}$ with $t_k >0$ satisfying \eqref{eq:recurrence_OPS}; 
\item[$(8)$] There exist constants $\sigma=(s_k)_{k\geq 0}$ and $\tau = (t_k)_{k\geq 1}$ with $t_k >0$ satisfying \eqref{recursive_matrix};
\item[$(9)$] There exists a positive-definite Riesz functional $\mathcal{L}_y$ satisfying \eqref{eq:Riesz_condition}.
\end{enumerate}
\end{theorem}
\begin{proof}
Note that $y$ is the Catalan-like numbers if and only if $\Delta_m (y) \neq 0$ for all $m\in \N$. By the definition of the positive definite sequence,  (1), (2), and (3) are equivalent.
By Theorem \ref{thm:widder1},  (2) and (4) are equivalent. By a similar proof of Theorem \ref{thm:main1},  it is easy to check the remainder.
\end{proof}

If a sequence is positive, then there are two disjoint classes as follows.
\begin{itemize}
\item[(i)] The sequence $y$ is \emph{positive definite} if the Hankel matrix $H_m(y)$ is positive definite for all $n\in \N$.
\item[(ii)] The sequence $y$ is \emph{positive semidefinite} if it is not positive definite.
\end{itemize}
Observe that if a sequence is positive semidefinite, then the Hankel matrix $H_m(y)$ is positive semidefinite for all $n \in \N$ and at least one of them is singular. In fact, once one of the finite Hankel matrix is singular, all the following ones are also singular. Based on this definition, the sets of positive definite sequences and positive semidefinite sequences are mutually disjoint.

% Theorem: Positive definite & positive semidefinite
\begin{theorem}[\cite{book:widder}, Theorem~12a]\label{thm:widder1}
A necessary and sufficient condition that there exists a nonnegative Borel measure $\mu$ with $\big| \supp(\mu) \big| = \infty$ (resp. $\big| \supp(\mu) \big| < \infty$) such that
\begin{equation}\label{eq:integral_mu}
y_{k}=\int_{\R} x^k \mathrm{d}\mu \quad \textup{for all } k\in \N
\end{equation}
is that the sequence $y$ is positive definite (resp. positive semidefinite).
\end{theorem}

When the Riesz functional is positive, there is an intimate relationship between the zeros of the corresponding orthogonal polynomials
and the support of representing measure.
We maintain the hypothesis that the Riesz functional $\mathcal{L}$ is positive.
For proof of each theorem and proposition, see \cite{Chihara78} and references therein.

We denote the zeros of  $P_n(x)$ by $x_{n_i}$ with 
$$ x_{n1} < x_{n2} < \cdots < x_{nn}.$$
We denote 
\begin{equation*}
U=\{\xi_i | i=1,2,3,\ldots\}, \quad V=\{\eta_j | j=1,2,3,\ldots\},
\end{equation*}
where 
\begin{equation*}
\xi_i = \lim_{n\rightarrow \infty} x_{ni}, \quad \eta_j = \lim_{n\rightarrow \infty} x_{n,n-j+1}  \quad i,j =1,2,3,\ldots.
\end{equation*}

Since $\xi_{i-1} \leq \xi_i < \eta_j \leq \eta_{j-1}$, we define 
\begin{align*}
\xi&=
\begin{cases}
 -\infty  \quad \text{if } \xi_i = -\infty \quad \text{for all }i=1,2,\ldots,  \\
  \lim_{i\rightarrow \infty} \xi_{i}   \quad \text{otherwise}.
\end{cases}\\
\eta&=
\begin{cases}
 +\infty  \quad \text{if } \eta_j = +\infty \quad \text{for all }j=1,2,\ldots,  \\
  \lim_{j\rightarrow \infty} \xi_{j}   \quad \text{otherwise}.
\end{cases}
\end{align*}
Taking $\xi_0 = -\infty$ and $\eta_0 = + \infty$, we have
\begin{equation*}
-\infty = \xi_0 \leq \xi_1 \leq \xi_2 \leq \cdots \leq \xi \leq \eta \leq \cdots \leq \eta_2 \leq \eta_1 \leq \eta_0=+\infty.
\end{equation*}

%\xi & =  \lim_{i\rightarrow \infty} \xi_{i}, \quad \eta =  \lim_{j\rightarrow \infty} \eta_{j}.

Note that when $\xi$ and $\eta$ are both finite, $\supp(\mu) = U \cup S \cup V$ with $S\subseteq [\xi, \eta]$. 
For more information, see \cite[p.63]{Chihara78}.

For each $n\in \N$ we denote 
$$\alpha_n(x) = \frac{t_{n+1}}{(s_n - x)(s_{n+1} -x)}.$$

%The Catalan-like numbers are closely related to continued fractions and orthogonal polynomials.  Let $y_n$
%be the Catalan-like numbers corresponding to $(\sigma,\tau)$. Then 
%\begin{equation*}
%\sum_{n\geq 0} y_n x^n = \frac{1}{1-s_0x-\frac{t_1x^2}{1-s_1x-\frac{t_2x^2}{1-s_2x-\cdots}}}.
%\end{equation*}

\begin{definition}
The closed interval, $[\xi_1, \eta_1]$, is called the \emph{true interval of orthogonality} of the OPS (of $\mathcal{L}$).
\end{definition}
The true interval of orthogonality is the smallest closed interval that contains all of the zeros of all $P_n(x)$. Obviously the true interval of orthogonality is the smallest closed interval that includes $\supp(\mu)$(see \cite[Theorem 3.2]{Chihara78}).

The theory of chain sequences can provide a relation between the true interval of orthogonality, $[\xi_1, \eta_1]$, and the sequences 
$\sigma=(s_k)_{k\geq 0}$ and $\tau = (t_k)_{k\geq 1}$.

\begin{definition}
A sequence  $y=(y_n)_{n \geq 0}$ is called a $\emph{chain sequence}$ if there exists a sequence  $(g_k)_{k\geq 0}$ such that
\begin{itemize}
\item[$(i)$] $0 \leq g_0 < 1$,
\item[$(ii)$] $0< g_{n+1} < 1$ for all $n\in \N$,
\item[$(iii)$] $y_{n} = (1-g_{n}) g_{n+1}$ for all $n\in \N$.
\end{itemize}
Such  $(g_k)_{k\geq 0}$ is called a \emph{parameter sequence} for $y$, and 
 $g_0$ is called an \emph{initial parameter}.
\end{definition}

For the proof of the following theorems, see \cite[Theorem 2.1--2.2, Chapter IV]{Chihara78}.
\begin{theorem}\label{thm:chain}
$\xi_1 \geq a$ if and only if
\begin{itemize}
\item[$(i)$] $s_n > a $ for all $n \in \N$,
\item[$(ii)$] $(\alpha_n(a))_{n\geq 0}$ is a chain sequence.
\end{itemize}
\end{theorem}

\begin{corollary}\label{cor:chain}
$\eta_1 \leq b$ if and only if
\begin{itemize}
\item[$(i)$] $t_n < b $ for all $n \in \N$,
\item[$(ii)$] $(\alpha_n(b))_{n\geq 0}$ is a chain sequence.
\end{itemize}
\end{corollary}
We will mainly use Theorem \ref{thm:chain} and Corollary \ref{cor:chain} to prove Theorem \ref{thm:support1}, which is one of our main results.

%\begin{theorem}\label{thm:bounded}
% $[\xi_1, \eta_1]$ is bounded if and only if both
%$(s_k)_{k\geq 0}$ and $(t_k)_{k\geq 1}$ are bounded.
%\end{theorem}

%\begin{theorem}
%Each Catalan-like numbers in Example \ref{ex:Catalan-like} are Hausdorff moment sequences. Furthermore, the support  
%\end{theorem}
%\begin{proof}
%By Theorem \ref{thm:bounded} it is easy to check.
%

%
%\begin{itemize}
%\item[$(i)$] For the Catalan numbers $C_n$  $\sigma  = (1,2,2,\ldots)$ and $\tau  = (1,1,1,\ldots)$.
%We find $b\in \R$ such that $(\alpha_n(b))_{n\geq 0}$ is a chain sequence, i.e.,
%\begin{equation*}
%\begin{cases}
%\alpha_0(b) = \dfrac{1}{(1- b)(2 -b)} = (1-g_0)g_1\\
%\alpha_n(b) = \dfrac{1}{(2- b)^2} = (1-g_n)g_{n+1} \quad \text{for all } n\in \N.
%\end{cases}
%\end{equation*}
%Assuming $g_{n+1}=g_{n+2}$ for all $n \in \N$, it holds that $b \geq 4$. Clearly, $t_n < b$ for all $n$ and $(\alpha_n(4))_{n\geq 0}$ is a chain sequence. 
%Indeed, the parameter sequence for  $(\alpha_n(4))_{n\geq 0}$ is 
%$\frac{2}{3},~ \frac{1}{2}, ~ \frac{1}{2}, ~ \frac{1}{2}, \ldots.$
%Thus, by Corollary \ref{cor:chain} $\eta_1 \leq 4$. Therefore, the support of measure for an integral representation of Catalan numbers is subset of $[0,4]$.
%Using simple linear transform, the support can be converted into subset of $[0,1]$.
%For more information about integral representations of the Catalan numbers, see the survey article\cite{math5030040} and therein.
%\end{itemize}
%\end{proof}

% Main results

\section{Catalan-like numbers with $y_n(p,s ; q,t)$}
As a main result, we given an affirmative answer to the question of \cite{LIANG2016484}, we will show that many well-known Catalan-like numbers are Hausdorff moment sequences.

We denote by $y(p,s ; q,t)$ the Catalan-like numbers corresponding to
$ \sigma = (p, s, s, \ldots)$ and $\tau = (q, t, t, \ldots)$.
Remark that many well-known Catalan-like numbers are of the such form.

\begin{theorem}\label{thm:support1}
%Let $p,s,q,t$ be nonnegative numbers.
If $p> s-2\sqrt{t}$ and $\max\{q,t\} < s+2\sqrt{t}$, then 
the Catalan-like numbers $y(p,s;q,t)$  form a Hausdorff moment sequence such that 
the support of representing measure is contained in $[s-2\sqrt{t}, s+2\sqrt{t}]$.
%Moreover, $\{\xi_i | i=1,2,3,\ldots\} =\{s-2\sqrt{t}\}$, $\{\eta_j | j=1,2,3,\ldots\}=\{s+2\sqrt{t}\}$, $\xi = s-2\sqrt{t}$, and $\eta = s+2\sqrt{t}$.
In fact,  $[s-2\sqrt{t}, s+2\sqrt{t}]$ is the smallest interval including the support.
Furthermore, if $y(p,s;q,t)$ additionally holds the condition, $s\geq 2\sqrt{t}$, then it is also Stieltjes moment sequence as well.
\end{theorem}
\begin{proof}
For the Catalan-like numbers $y_n$ corresponding to  $\sigma  = (p,s,s,\ldots)$ and $\tau  = (q,t,t,\ldots)$.
We find $x\in \R$ such that $(\alpha_n(x))_{n\geq 0}$ is a chain sequence, i.e.,
\begin{equation*}
\alpha_n(x)=
\begin{cases}
 \dfrac{q}{(p- x)(s -x)} = (1-g_0)g_1  \quad \text{if } n=0 \\
  \dfrac{t}{(s- x)^2} = (1-g_n)g_{n+1}  \quad \text{otherwise }n \geq 1.
\end{cases}
\end{equation*}
Assume that $g_n$ is convergent to the limit, called $g$.
%Let $\lim_{n \rightarrow \infty} g_n =g$.
Then we have  
$$ \lim_{n \rightarrow \infty} \alpha_n(x) = \frac{t}{(s- x)^2} = (1-g)g. $$
Since $0 \leq g \leq 1$, it holds that $x\leq s-2\sqrt{t}$ or $x\geq s+2\sqrt{t}$.
By Theorem \ref{thm:chain} and Corollary \ref{cor:chain}, it follows that 
 $[\xi_1, \eta_1] \subseteq [s-2\sqrt{t}, s+2\sqrt{t}]$.

Using the fact in \cite[p.121]{Chihara78}, we have $\xi = s-2\sqrt{t}$ and $\eta = s+2\sqrt{t}$.
Thus, $s-2\sqrt{t} \leq \xi_i \leq \xi =s-2\sqrt{t}$ and $s+2\sqrt{t} = \eta \leq \eta_j \leq s+2\sqrt{t}$.
\end{proof}
It is easy to check that
$$ 1-\frac{q}{\sqrt{t}(p-s+2\sqrt{t})},~ \frac{1}{4},~ \frac{1}{4},~ \frac{1}{4},~ \ldots$$
is a parameter sequence for $\alpha_n(s-2\sqrt{t})$.

%In fact, it holds that
%\begin{equation*}
%-\infty = \xi_0 < \xi_1 = \xi_2 = \cdots = \xi = s-2\sqrt{t} \leq s+2\sqrt{t} = \eta= \cdots = \eta_2 = \eta_1 < \eta_0=+\infty,
%\end{equation*}
%where $\xi_i, \eta_i, \xi, \eta$ are defined in 

Theorem 3.13 in \cite{book:momentproblem} states that a sequence $y$ is an $[a,b]$-moment sequence if and only if
$H_m(y)$ and $(a+b)H_m(Ey)-H_m(E(Ey))-abH_m(y)$ are positive semidefinite for all $m\in \N$.
Here $Ey$ denote the shifted sequence given by $(Ey)_n=(y_{n+1})$, $n\in \N$.

\begin{corollary}\label{cor:necessary_condition}
\begin{itemize}
\item[(1)] When $y=y_n(p,2;q,1)$ satisfy $p>0$ and $0<q<4$, then it has a representing measure with the support contained in  $[0,4]$, which is equivalent to $H_m(y)\geq 0$ and $H_m(4Ey-E(Ey))\geq 0$ for all $m\in \N$.

\item[(2)] When $y=y_n(p,1;q,1)$ satisy $p>-1$ and $0<q<3$, then it has a representing measure with the support contained in  $[-1,3]$, which is equivalent to $H_m(y)\geq 0$ and $H_m(3y+2Ey-E(Ey))\geq 0$ for all $m\in \N$.

\item[(3)] When $y=y_n(p,3;q,2)$ satisy $p>3-2\sqrt{2}$ and $0<q<3+2\sqrt{2}$, then it has a representing measure with the support contained in  $[3-2\sqrt{2},3+2\sqrt{2}]$, which is equivalent to $H_m(y)\geq 0$ and $H_m(-y+6Ey-E(Ey))\geq 0$ for all $m\in \N$.

\item[(4)] When $y=y_n(p,3;q,1)$ satisy $p>1$ and $0<q<5$, then it has a representing measure with the support contained in  $[1,5]$, which is equivalent to $H_m(y)\geq 0$ and $H_m(-5y+6Ey-E(Ey))\geq 0$ for all $m\in \N$.

\end{itemize}
\end{corollary}

Surprisingly, our results for support intervals of the Catalan-like numbers perfectly match the existing results for integral representations of well-known counting numbers.
\begin{example}
The Catalan-liken numbers in Example \ref{ex:Catalan-like} are Hausdorff moment sequences.
\begin{itemize}
\item[$(i)$] The Catalan numbers $C_n$,  the shifted Catalan numbers $C_{n+1}$,
the central binomial coefficients $\binom{2n}{n}$, and the Fine numbers $F_n$ 
have representing measures with compact support contained in  $[0,4]$.
For example, the Catalan numbers $C_n$  and the central binomial coefficients $\binom{2n}{n}$ are uniquely represented by 
\begin{equation*}
C_n = \int_0^4  x^n \Bigg( \frac{1}{2\pi}\sqrt{\frac{4-x}{x}}  \Bigg) dx
\quad \text{and} \quad \binom{2n}{n} = \int_0^4  x^n  \Bigg(\frac{1}{\pi \sqrt{x(4-x)}} \Bigg) dx,
\end{equation*}
respectively.
For more information about integral representations of the Catalan numbers, see \cite{math5030040} and references therein.
%And it holds that 
%\begin{equation*}
%\begin{bmatrix}
%4C_{1}-C_{2} & 4C_{2}-C_{3} & \cdots & 4C_{m+1}-C_{m+2}  \\
%4C_{2}-C_{3}  & 4C_{3}-C_{4}  & \cdots & 4C_{m+2}-C_{m+3}  \\
%\vdots & \vdots & \ddots & \vdots \\
%4C_{m+1}-C_{m+2}  & 4C_{m+2}-C_{m+3}  & \cdots & 4C_{2m+1}-C_{2m+2}  \\
%\end{bmatrix}.
%\end{equation*}
%is positive semidefinite for all $m\in \N$.

\item[$(ii)$] The Motzkin numbers $M_n$,
 the central trinomial coefficients $T_n$, and
the Riordan numbers $R_n$ have representing measures with compact support contained in $[-1,3]$.
For example, the Motzkin numbers $M_n$ and the central trinomial coefficients $T_n$ are represented by 
\begin{equation*}
M_n = \int_{-1}^3   x^n  \Bigg(\frac{1}{2\pi} \sqrt{(3-x)(1+x)} \Bigg) dx
\quad \text{and} \quad
%And gives the integral representation for the central trinomial coefficients $T_n$ as follows.
T_n = \int_{-1}^3   x^n  \Bigg(\frac{1}{\pi \sqrt{(3-x)(1+x)}} \Bigg) dx,
\end{equation*}
respectively (See \cite{PRB11}).
%And it holds that 
%\begin{equation*}
%\begin{bmatrix}
%2M_{1}-M_{2}+3M_{0} & 2M_{2}-M_{3}+3M_{1} & \cdots & 2M_{m+1}-M_{m+2}+3M_{m}  \\
%2M_{2}-M_{3}+3M_{1}  & 2M_{3}-M_{4}+3M_{2}  & \cdots & 2M_{m+2}-M_{m+3}+3M_{m+1}  \\
%\vdots & \vdots & \ddots & \vdots \\
%2M_{m+1}-M_{m+2}+3M_{m}  & 2M_{m+2}-M_{m+3}+3M_{m+1}  & \cdots & 2M_{2m+1}-M_{2m+2}+3M_{2m}  \\
%\end{bmatrix}.
%\end{equation*}
%is positive semidefinite for all $m\in \N$.

\item[$(iii)$] The central Delannoy numbers $D_n$,
the large Schr\'{o}der numbers $r_n$, and
the little Schr\'{o}der numbers $S_n$
have representing measures with compact support contained in $[3-2\sqrt{2},3+2\sqrt{2}]$.
 \begin{equation*}
D_n = \int_{3-2\sqrt{2}}^{3+2\sqrt{2}}   x^n  \Bigg(\frac{1}{\pi} \frac{1}{\sqrt{\big(3+2\sqrt{2}-x \big) \big(x-3+2\sqrt{2} \big)}} \Bigg) dx.
\end{equation*}
This integral representations can be obtained from the integral form \cite[Theorem 1.3]{QI2018101} 
% And it holds that 
%\begin{equation*}
%\begin{bmatrix}
%6D_{1}-D_{2}-D_{0} & 6D_{2}-D_{3}-D_{1} & \cdots & 6D_{m+1}-D_{m+2}-D_{m}  \\
%6D_{2}-D_{3}-D_{1} & 6D_{3}-D_{4}-D_{2}  & \cdots & 6D_{m+2}-D_{m+3}-D_{m+1}  \\
%\vdots & \vdots & \ddots & \vdots \\
%6D_{m+1}-D_{m+2}-D_{m} & 6D_{m+2}-D_{m+3}-D_{m+1}  & \cdots & 6D_{2m+1}-D_{2m+2}-D_{2m}  \\
%\end{bmatrix}.
%\end{equation*}
%is positive semidefinite for all $m\in \N$.

\item[$(iv)$] The (restricted) hexagonal numbers $h_n$ has representing measures with compact support contained in $[1,5]$.
% And it holds that 
%\begin{equation*}
%\begin{bmatrix}
%6h_{1}-h_{2}-5h_{0} & 6h_{2}-h_{3}-5h_{1} & \cdots & 6h_{m+1}-h_{m+2}-5h_{m}  \\
%6h_{2}-h_{3}-5h_{1} & 6h_{3}-h_{4}-5h_{2}  & \cdots & 6h_{m+2}-h_{m+3}-5h_{m+1}  \\
%\vdots & \vdots & \ddots & \vdots \\
%6h_{m+1}-h_{m+2}-5h_{m} & 6h_{m+2}-h_{m+3}-5h_{m+1}  & \cdots & 6h_{2m+1}-h_{2m+2}-5h_{2m}  \\
%\end{bmatrix}.
%\end{equation*}
%is positive semidefinite for all $m\in \N$.
\end{itemize}
\end{example}

Many well-known combinatorial numbers can be expressed as integrals with $[a,b]$-representing measures.
Remark that every Hausdorff moment sequence is determinate. Thus the Catalan-liken numbers in Example \ref{ex:Catalan-like} are determinate, which means that for each sequence there is the unique measure, respectively. 
In \cite{BergSzwarc} Berg and Szwarc prove that if a sequence $y$ satisfies $\Delta_m(y)>0$ for $m<s$ while $\Delta_m(y)=0$ for $m\geq s$, then all Hankel matrices are positive semidefinite, and in particular, $y$ is a Hamburger moment sequence  with a discrete measure $\mu$ such that $\big| \supp(\mu) \big| = s$.
Remark that by Theorem \ref{thm:main1}, it holds $\Delta_m (y) \neq 0$, $m\geq 0$ for any Catalan-like numbers. Thus the supports of representing measures are not finite.
Thus, the representing measures for Catalan-like numbers cannot be expressed as simple finite discrete measures. 
Although we do not know how to obtain the representing measure in general, Mnatsakanov\cite{MNATSAKANOV20081612} provided an approximation of the measure with a given Hausdorff moment sequence. Remark that to find such approximations is required to check if a given sequence is a Hausdorff moment sequence and to find the interval including the support.
Our result can be helpful to do it.

\section{New sequences generated by Catalan-like numbers}

\subsection{Subsequences of Catalan-like numbers}

\begin{example}\label{ex:catalan1}
Consider the following sequence  $y=(y_n)_{n\geq 0}$.
\begin{equation*}
1,~ 0,~ 1,~ 0,~ 2,~ 0,~ 5,~ 0,~ 14,~ 0,~ 42,~ 0,~ 132,~ 0,~ 429,~ 0 , \ldots.
\end{equation*}
It is easy to check that $(y_n)_{n\geq 0}$ is the 
Catalan-like numbers corresponding to $(\sigma, \tau)$ with $\sigma  = (0,0,0,\ldots)$ and $\tau  = (1,1,1,\ldots)$. 
Note that since $\det(H_m(y))>0$ for all $m\in \N$ and $\det(\widetilde{H}_2(y))=-1$, $(y_n)_{n\geq 0}$ is a Hamburger moment sequence, but it is not a Stieltjes moment sequence.
 However,  its subsequence, $(y_{2n})_{n\geq 0}$, is a Stieltjes moment sequence which is the Catalan numbers. 
 \end{example}
 This example motivates us to study subsequences of Catalan-like numbers.
In \cite{ACJ17} 
for a given Hamburger moment sequence $(y_n)_{n\geq 0}$, it is shown that its subsequence $(y_{n_k})_{k \geq 0}$, given by $n_k=dk+\ell$ ($d\in \N,~ \ell\in 2\N$) for all $k\in\N$ is always a Hamburger moment sequence. Also, a relationship between Hamburger moment sequence and its subsequences via Cauchy transform is provided.  For multisequences, see \cite{SHF17}.
In this paper we consider subsequences of Catalan-like numbers which admit $K$-measures when $K=\R$, $K=[0,\infty)$, or $K=[0,1]$.

By a subsequence of a sequence $(y_n)_{n\geq 0}$, we shall mean a
sequence of the form $(y_{n_k})_{k \geq 0}$, where each $n_k\in \N$ and $n_0 < n_1 < \cdots.$

Since almost all of the well-known Catalan-like numbers are positive sequences, 
we may find some connections between their subsequences and 
consider subsequences of the classical $K$-moment sequences for $K\subseteq \R$. 
For given $d,\ell\in \N$ consider the subsequence $(y_{n_k})_{k\geq0}$, defined by 
\begin{equation}\label{eq:subsequence1}
n_k = dk+\ell \quad \text{for all }k \in \N.
\end{equation}

\begin{theorem}\label{thm:submoments1Dmain}
Let $K\subseteq \R$ be a closed set.
Let $d\in \N$ and let
\begin{equation}\label{eq:ell}
\begin{cases}
\ell \in \N & \text{if } K \subseteq [0,\infty),\\
\ell \in 2 \N & \text{otherwise.}
\end{cases}
\end{equation}
If $y=(y_n)_{n \geq 0}$ is a $K$-moment sequence, then subsequence, 
$\tilde{y}=(\tilde{y}_{k})_{k \geq 0}$, defined by $\tilde{y}_k = y_{dk+\ell}$ for all $k\in\N$,  is a $\widetilde{K}$-moment sequence, where $\widetilde{K} = \{ x^d | x\in K \}$. \end{theorem}
\begin{proof}
Note that $\widetilde{K}$ is a closed set.
%Let $K=[0,1]$ or $K=[0,\infty)$. 
Since $y$ admits a $K$-representing measure, the Riesz functional $\mathcal{L}_y$ is $K$-positive. We will show that $\mathcal{L}_{\tilde{y}}(p) \geq 0$ for all 
$p\in \R[\bm{x}]$ with $p|_{\widetilde{K}} \geq 0$.
Let $p(x)=\sum p_k x^k \in \R[x]$ such that $p |_{\widetilde{K}} \geq 0$. 
%and $d,\ell \in \N$.
Since $q(x):=p(x^d)x^{\ell} \geq 0$ on $K$, it follows that
$$\mathcal{L}_{\tilde{y}}(p) = \sum p_k \tilde{y}_{k} =  \sum p_k y_{dk+\ell} = \mathcal{L}_{y}(q) \geq 0.$$
Thus, $\mathcal{L}_{\tilde{y}}$ is $K$-positive. Therefore, by Theorem \ref{thm:RH_full} $\tilde{y}$ admits a $\widetilde{K}$-representing measure.
\end{proof}

\begin{corollary}\label{thm:submoments1D}
If $y=(y_n)_{n \geq 0}$ is a $(A)$-moment sequence, then subsequence, 
$\tilde{y}=(\tilde{y}_{n})_{n \geq 0}$, is a $(B)$-moment sequence with respect to even and odd numbers of $d$ and $\ell$ as follows:
\begin{center}
\begin{tabular}{ |c|c|c|c||c|c|c|c| } 
\hline
$d$ & $\ell$ & $(A)$ & $(B)$ & $d$ & $\ell$ & $(A)$ & $(B)$\\
\hline
 &  & $(-\infty,\infty)$ & N/A  & & & $(-\infty,\infty)$  & N/A  \\ \cline{3-4} \cline{7-8}
 & odd & $[0,\infty)$ &  $[0,\infty)$ &  & odd &  $[0,\infty)$ &  $[0,\infty)$\\ \cline{3-4}\cline{7-8}
odd &  & $[0,1]$  &  $[0,1]$ & even &  & $[0,1]$ &  $[0,1]$\\ \cline{2-4}\cline{6-8}
 &  & $(-\infty,\infty)$  & $(-\infty,\infty)$   &  & & $(-\infty,\infty)$  &  $[0,\infty)$\\ \cline{3-4}\cline{7-8}
& even &  $[0,\infty)$ & $[0,\infty)$ &  & even &  $[0,\infty)$ &  $[0,\infty)$\\ \cline{3-4}\cline{7-8}
&  &  $[0,1]$ & $[0,1]$ &  & & $[0,1]$ &  $[0,1]$ \\ \cline{3-4}\cline{7-8}
\hline
\end{tabular}
\end{center}

{\bf Summary}
(i) If $y$ is a Hausdorff (resp. Stieltjes)  moment sequence, then
$\tilde{y}$ is a Hausdorff moment sequence (resp. Stieltjes) for all $d,\ell\in \N$.
(ii) If $y$ is Hamburger moment sequence, then
$\tilde{y}$ is a Hamburger moment sequence for all $d\in 2\N+1$ and $\ell\in 2\N$.
(iii) If $y$ is Hamburger moment sequence, then
$\tilde{y}$ is a Stieltjes moment sequence for all $d\in 2\N$ and $\ell\in 2\N$.
(iv) If $y$ is Hamburger moment sequence, then
$\tilde{y}$ may or may not be a Hamburger moment sequence for all $d\in \N$ and $\ell\in 2\N+1$.
\end{corollary}
%{\bf Summary}
%\begin{itemize}
%    \item[(i)] If $y$ is a Hausdorff (resp. Stieltjes)  moment sequence, then
%$\tilde{y}$ is a Hausdorff moment sequence (resp. Stieltjes) for all $d,\ell\in \N$.
%
%
%\item[(ii)] If $y$ is Hamburger moment sequence, then
%$\tilde{y}$ is a Hamburger moment sequence for all $d\in 2\N+1$ and $\ell\in 2\N$.
%
%\item[(iii)] If $y$ is Hamburger moment sequence, then
%$\tilde{y}$ is a Stieltjes moment sequence for all $d\in 2\N$ and $\ell\in 2\N$.
%
%\item[(iv)] If $y$ is Hamburger moment sequence, then
%$\tilde{y}$ may or may not be a Hamburger moment sequence for all $d\in \N$ and $\ell\in 2\N+1$.
%\end{itemize}

\begin{proof}
By Theorem \ref{thm:submoments1Dmain} it is easy to check that (i), (ii), (iii) are true.
(iv) If $y$ is a Hamburger moment sequence which is a Stiltjes moment sequence, then it is trivial that $\tilde{y}$ is a Hamburger moment sequence.
Consider a Hamburger moment sequence which is not a Stieltjes moment sequence.
Let $\mu=\frac{1}{2}\delta_{-1} + \frac{1}{2}\delta_{2}$ and let 
\begin{equation*}
%y_k = \frac{1}{2} ((-1)^k+2^k) = \int_\R x^k d\mu(x) \quad \text{for all } k\in \N,
y_k  = \int_\R x^k d\mu \quad \text{for all } k\in \N.
\end{equation*}
Since $y$ is a Hamburger moment sequence, it holds $H_m(y) \geq 0$ for all $m\in \N$.
However, its subsequence $\tilde{y}$ with $d=1$ and $\ell=1$ does not admit any $K$-representing measure with $K=\R$. Indeed,
$\det([\tilde{y}_{i+j}]_{0 \leq i,j \leq 1}) = -\frac{9}{2} <0$.

%\begin{itemize}
%\item[(i)] Let $K=[0,1]$ or $K=[0,\infty)$. Since $y$ is a $K$-moment, the Riesz functional $\mathcal{L}_y$ is $K$-positive.
%Let $p(x)=\sum p_k x^k \in \R[x]$ such that $p|_{K} \geq 0$ and $d,\ell \in \N$.
%Since $q(x):=p(x^d)x^{\ell} \geq 0$ on $K$, it follows that
%$$\mathcal{L}_{\tilde{y}}(p) = \sum p_k \tilde{y}_{k} =  \sum p_k y_{dk+\ell} = L_{y}(q) \geq 0.$$
%Thus, $\mathcal{L}_{\tilde{y}}$ is $K$-positive. Therefore, by Theorem \ref{thm:RH_full} $\tilde{y}$ is a $K$-moment sequence.
%
%\item[(ii)] Note that if  $p\in \R[x]$ and $p|_K \geq 0$ for $K=\R$, then $q(x)=p(x^d)x^{\ell} \geq 0$ on $K=\R$ for all $d\in 2\N+1$ and $\ell\in 2\N$.
%
%\item[(iii)] Let $K=\R$ and $\tilde{K}=[0,\infty)$. Since $y$ is a $K$-moment, the Riesz functional $\mathcal{L}_y$ is $K$-positive.
%Let $p(x)=\sum p_k x^k \in \R[x]$ such that $p|_{\tilde{K}} \geq 0$ and $d,\ell \in 2\N$.
%Since $q(x):=p(x^d)x^{\ell} \geq 0$ on $K$, it follows that
%$$\mathcal{L}_{\tilde{y}}(p) = \sum p_k \tilde{y}_{k} =  \sum p_k y_{dk+\ell} = \mathcal{L}_{y}(q) \geq 0.$$
%Thus, $\mathcal{L}_{\tilde{y}}$ is $\tilde{K}$-positive. Therefore, by Theorem \ref{thm:RH_full} $\tilde{y}$ is a $\tilde{K}$-moment sequence.
%\end{itemize}
\end{proof}

Remark that it is trivial that if $(y_n)_{n\geq 0}$ is a $[0,1]$-moment sequence, then it is a Stieltjes moment sequence and a  Hamburger moment sequence, and similarly if $(y_n)_{n\geq 0}$ is a Stieltjes moment sequence, then it is a Hamburger moment sequence.

\begin{theorem}\label{pattern1}
A necessary and sufficient condition for $n_k$ that the subsequence $(y_{n_k})_{n\geq 0}$ is Stieltjes moment sequences for all Stieltjes moment sequences $(y_n)_{n\geq 0}$ is
$$n_k=dk+\ell \quad \text{for all } k \in \N,$$
where $d,\ell\in \N$.
\end{theorem}
\begin{proof}
(Sufficient condition:) Let $y=(y_k)_{k\in\N}$ be a Stieltjes moment sequence. Then by Corollary \ref{thm:submoments1D}, it follows that $(y_{n_k})_{k\geq 0}$ is also Stieltjes moment sequence. 

(Necessary condition:) Suppose that there exists a $n_k$ such that $n_k$ is not of the above form, but  $(y_{n_k})_{k\geq 0}$ is Stieltjes moment sequence. Then there exist three consecutive terms $n_s,n_{s+1},n_{s+2} \in \N$ such that 
$n_{s+1} - n_s \neq n_{s+2}- n_{s+1}$. At first we assume that $n_{s+1} - n_s < n_{s+2}- n_{s+1}$.
Set $d= n_{s+1} - n_s$ and $a=n_s$. Then, $n_{s+1} = a+d$, and  $n_{s+2} = a+2d+e$ for $e \geq 1$.

Consider the Stieltjes moment sequence whose representing measure is $\mu=\delta_{\epsilon}$ with $\epsilon=1/2$.
Then, it follows that 
\begin{equation}\label{2by2submatrix}
    \begin{pmatrix}
    y_{n_s} &  y_{n_{s+1}}\\
     y_{n_{s+1}} &  y_{n_{s+2}}
    \end{pmatrix}
    =
     \begin{pmatrix}
    y_{a} &  y_{a+d}\\
     y_{a+d} &  y_{a+2d+e}
    \end{pmatrix}
    =
         \begin{pmatrix}
    \epsilon^a &  \epsilon^{a+d}\\
     \epsilon^{a+d} &  \epsilon^{a+2d+e}
    \end{pmatrix}
\end{equation}
is not positive definite, which is a contradiction. 
Note that unlike Hamburger moment sequence, 
Stieltjes moment sequence have the principal submatrix which is of the form (\ref{2by2submatrix}),
regardless of even and odd-ness of $a$.

When $n_{s+1} - n_s > n_{s+2}- n_{s+1}$, we can show a contradiction by a similar argument (use $\mu=\delta_{\epsilon}$ with $\epsilon=2$).
\end{proof}

Remark that although $y=(y_n)_{n\geq 0}$ is the Catalan-like numbers, its subsequence $(y_{n_k})_{k\geq0}$, defined by $n_k = dk+\ell$ for all $ k \in \N$, is possibly not 
the Catalan-like numbers. For example, the sequence in Example \ref{ex:catalan1} is the Catalan-like numbers, but its subsequence $(\tilde{y}_{k})_{k\geq0}$, defined by $\tilde{y}_k = y_{2k+1}$, is not the Catalan-like numbers, since $\Delta_m(\tilde{y}) \neq 0$ for all $m\in \N$.

A sequence $y=(y_n)_{n \geq 0}$ is called \emph{Stieltjes (resp. Hausdorff) Catalan-like numbers}
if it is a Stieltjes (resp. Hausdorff) moment sequence and the Catalan-like numbers.

\begin{theorem}
If $y=(y_n)_{n\geq 0}$ is Stieltjes Catalan-like numbers corresponding to $(\sigma,\tau)$, then 
the subsequence $\tilde{y}=(y_{n_k})_{k\geq0}$, defined by $n_k = dk+\ell$ for all $ k \in \N$, is
Stieltjes Catalan-like numbers corresponding to $\tilde{\sigma} = (\tilde{s}_k)_{k\geq 0}$ and $\tilde{\tau} = (\tilde{t}_k)_{k\geq 1}$, where 
\begin{equation}
\tilde{s}_k = \frac{\mathcal{L}_{\tilde{y}}[x\tilde{P}_k^2(x)]}{\mathcal{L}_{\tilde{y}}[ \tilde{P}_k^2(x)]} \quad \text{and} \quad \tilde{t}_k = \frac{\mathcal{L}_{\tilde{y}}[\tilde{P}_k^2(x)]}{\mathcal{L}_{\tilde{y}}[\tilde{P}_{k-1}^2(x)]}
% = \frac{\Delta_{k-2} \Delta_{k}}{\Delta_{k-1}^2}.
\end{equation}
such that
\begin{equation*}%\label{Hankelmatrix}
\tilde{P}_n(x)=
\frac{1}{\Delta_{n-1}(\tilde{y})}
\det
\begin{bmatrix}
y_{\ell} & y_{d+\ell} & \cdots & y_{dn + \ell} \\
y_{d+\ell} & y_{2} & \cdots & y_{d(n+1) + \ell} \\
\vdots & \vdots & \ddots & \vdots \\
y_{d(n-1)+\ell} & y_{dn+ \ell} & \cdots & y_{d(2n-1) +\ell} \\
1 & x & \cdots & x^n \\
\end{bmatrix}.
\end{equation*}
\end{theorem}
\begin{proof}
By Corallary \ref{thm:submoments1D} (1), the subsequence $\tilde{y}$ is a Stieltjes moment sequences. Since $y$ is Stieltjes Catalan-like numbers, it is a positive sequence and Catalan-like numbers. By Theorem \ref{thm:main2} (3), it holds that $H_m(y)>0$ for all $m\in \N$. Let $\lambda_m$ and $\tilde{\lambda}_m$ be the smallest eigenvalue of $H_m(y)$ and $\tilde{H}_m(\tilde{y})$. Since $\lambda_m>0$ for all $m\in \N$, by Cauchy interlacing theorem it holds $\tilde{\lambda}_m>0$ for all $m\in \N$.
By Theorem  \ref{thm:OPS_form} and \ref{thm:recurrence_OPS}, it is easy to find the explicit forms of $\tilde{\sigma}$ and $\tilde{\tau}$.
\end{proof}

%
%Even for $K \subseteq \R^n$ being not of the form $K=\prod_{i=1}^n K_i$, a similar statement holds as follows.
%\begin{proposition}
%Let $K\subseteq [0,\infty)^n$ be a closed set.
%If $(y_{\bm{\alpha}})$ admits a $K$-representing measure, then
%its subsequence $(\tilde{y}_{\bm{\alpha}})$ of the form (\ref{eq:subsequence})
%admits a  $K$-representing measure.
%\end{proposition}
%\begin{proof}
%Note that if $p\in \R[\bm{x}]$ and $p|_K \geq 0$ for $K=[0,\infty)^n$, then $q(\bm{x})=p(\bm{x}^{\bm{d}})x^{\bm{\ell}} \geq 0$ on $K$ for all $\bm{d},{\bm{\ell}} \in \N^n$.
%\end{proof}

%
%\begin{theorem}
%If $s_0\geq 1$ and $s_k \geq t_k +1$ for $k\geq 1$, then the Catalan-like numbers corresponding to $(\sigma, \tau)$ are s Stieltjes moment sequences.
%\end{theorem}
%
%\begin{remark}
%\begin{itemize}
%\item[(1)]
%% If the Catalan-like numbers, $y=(y_n)_{n\geq 0}$, admits a $K$-moment sequence with $K\subseteq \R$, then 
%%there exists $K$-representing measure for $y$ whose support is not finite. 
%
%\item[(2)] It is questionable that if the well-known combinatorial sequences in Example \ref{ex:Catalan-like} are Hausdorff moment sequences.
%\end{itemize}
%\end{remark}

%\section*{Remark}

\subsection{linear combinations of Catalan-like numbers}

Now we consider new sequences which are linear combinations of consecutive Catalan-like numbers.

Let $g(x)=\sum_{k=0}^n g_{n,k} x^k  \in \R[x]$.
For a sequence  $y=(y_k)_{k \geq 0}$, we define a new sequence $\mathcal{T}_g(y) = (\mathcal{T}_g(y)_k)_{k \geq 0}$ by
$$
 \mathcal{T}_g(y)_k= \sum_{k=0}^n g_{n,k} y_{n+k}.
 $$
\begin{theorem}
Let  $g(x)=\sum_{k=0}^n g_{n,k} x^k  \in \R[x]$ such that $g|_{[a,b]} \geq 0$.
If $y$ is an $[a,b]$-moment sequence,  
so is $\mathcal{T}_g(y)$. 
The measure for $\mathcal{T}_g(y)$ is $g\cdot d\mu$%$\tilde{\mu}(x)=g(x)\mu(x)$
, where $\mu$ is a representing measure for $y$.
\end{theorem}
\begin{proof}
Since  $y$ is an $[a,b]$-moment sequence, 
there exists a nonnegative measure, $\mu$,
such that
\begin{equation*}
y_{n} = \int_{a}^{b} x^{n} d\mu, \quad \text{for all } 
x \in \N.
\end{equation*}
%Set $z_k =  \int_{a}^{b} g(x) x^{n} d\mu$. 
Then it follows that
\begin{align*}
 \mathcal{T}_g(y)_k = \sum_{k=0}^n  g_{n,k} y_{n+k} =  \int_{a}^{b}  x^{k} g(x)    d\mu =  \int_{a}^{b}  x^{k}   d\tilde{\mu} ,
\end{align*} 
where $\tilde{\mu} = gd\mu$.
Since  $\mathcal{T}_g(y)$ is an $[a,b]$-moment sequence, it is determinate. Thus it has the unique representing measure.
\end{proof}
To see some results about the measure of subsequence, see \cite{ACJ17}.

\begin{example}
Let $y$ be an $[a,b]$-moment sequence. Then the following 
are also $[a,b]$-moment sequences.
\begin{itemize}
\item[(i)] $\mathcal{T}_g(y) = (\alpha y_k + \beta y_{k+1} )_{k\geq 0}$ 
with $g(x)= \alpha + \beta x$ such that $g|_{[a,b]} \geq 0$. 
\item[(ii)] $\mathcal{T}_g(y) = (-ab y_{k}  +(a+b) y_{k+1} - y_{k+2} )_{k\geq 0}$ 
with $g(x)= -(x-a)(x-b)$.
\item[(iii)] $\mathcal{T}_g(y) = (a^2b y_{k}  - (a^2+2ab) y_{k+1}  + (2a+b) y_{k+2}  -y_{k+3} )_{k\geq 0}$ 
with $g(x)= -(x-a)^2(x-b)$.
\end{itemize}
\end{example}

\begin{example}
Consider new sequences of the Catalan numbers as follows. 
\begin{itemize}
\item[$(i)$] ({\bf{Translation}}) The subsequence of Catalan numbers, $(C_{n+\ell})_{n\geq 0}$, is uniquely represented by 
\begin{equation*}
C_{n+\ell} = \int_0^4  x^n \Bigg( \frac{x^{\ell}}{2\pi}\sqrt{\frac{4-x}{x}}  \Bigg) dx
\end{equation*}

\item[$(ii)$] ({\bf{Moving with $d$ steps}})
The subsequence of Catalan numbers, $(C_{dn})_{n\geq 0}$, is uniquely represented by 
\begin{equation*}
C_{dn} = \int_0^{4^d}  x^n \Bigg( \frac{\sqrt[d]{x^{1-d}}}{2d\pi}  \sqrt{\frac{4-\sqrt[d]{x}}{ \sqrt[d]{x}}}  \Bigg) dx.
\end{equation*}

\item[$(iii)$] ({ \bf{Linear combinations} }) 
The new sequence $\tilde{C}:=(\tilde{C}_n)_{n\geq 0}=(4C_{n+1} - C_{n+2} )_{n\geq 0}$ 
is uniquely represented by 
\begin{equation*}
\tilde{C}_n = \int_0^4  x^n \bigg( \sqrt{x(4-x)^3} \bigg) dx.
\end{equation*}
\end{itemize}
%Note that $H_m(C_{dn+\ell})$ and $H_m(\mathcal{T}_g(y))$ are positive semidefinite for all $m\in \N$.
Note that they are all $[0,4]$-Hausdorff moment sequences. Thus they have unique representing measures, respectively.
\end{example}

\section{Remarks}
\begin{itemize}
\item[(i)]
Let $g(x)=\sum_{k=0}^n g_{n,k} x^k  \in \R[x]$.
If $g|_{[a,b]} \geq 0$, then $\mathcal{T}_g(y)$ is an $[a,b]$-moment sequence for any $[a,b]$-moment sequence $y$. In other words,
For a sequence  $y=(y_k)_{k \geq 0}$, the linear combinations of
consecutive Catalan-like numbers 
  $$
 \mathcal{T}_g(y)_k= \sum_{k=0}^n g_{n,k} y_{n+k}.
 $$
is  an $[a,b]$-moment sequence for any $[a,b]$-moment sequence $y$. 
It is questionable what is
a necessary and sufficient condition for a function $g\in \R[x]$ that $\mathcal{T}_g(y)$ is an $[a,b]$-moment sequence for all $[a,b]$-moment sequence $y$.

\item[(ii)] In this article we do not consider Catalan-like numbers which is not of the form $y(p,s ; q,t)$.
Some well-known counting numbers do not have such form. For instance,
the Bell numbers $B_n$ has $\sigma = \tau = (1,2,3,4,\ldots)$.
It was shown that the $B_n$ is Stieltjes moment sequence. However, we do not know whether it is a Hausdorff moment sequence.
\end{itemize}

\section{Acknowledgement}
This paper was completed at Institute of Mathematics, Academia Sinica, Taiwan during the visit of the first and the third authors. We would like to thank the institute for their hospitality and for research support. Without the support this work would not have been completed well.
The third-named author was supported by Basic Science Research Program through the National Research Foundation of Korea (NRF) funded by the Ministry of Education (2016R1A6A3A11932349).
The second-named author was supported partially by NSC under the Grant No.~107-2115-M-001-009-MY3.

\section*{References}
\bibliographystyle{plain}
\bibliography{references}

\end{document}